\documentclass[11pt]{article}

\usepackage{amsfonts}
\usepackage{amsthm,amsmath,amssymb,amsfonts}
\newtheorem{theorem}{{\bf Theorem}}[section]

\newtheorem{lemma}{{\bf Lemma}}[section]

\newtheorem{assumption}{{\bf Assumption}}[section]
\begin{document}
\title{\LARGE \bf Inverse Problem for a Curved Quantum Guide }
\date{}
\author{Laure Cardoulis\thanks{ Universit\'e de Toulouse, UT1 Ceremath, 21 All\'ees de Brienne, 
F-31042 Toulouse cedex, France, Institut de Math\'ematiques de Toulouse UMR 5219 ; laure.cardoulis@univ-tlse1.fr}
\and
Michel Cristofol\thanks{ Universit\'e d'Aix-Marseille, LATP, UMR 7353, 39, 
rue Joliot Curie, 13453 Marseille Cedex 13, France ; Michel.Cristofol@cmi.univ-mrs.fr}}

\maketitle
%
\numberwithin{equation}{section}

\abstract{In this paper, we consider the Dirichlet Laplacian operator $-\Delta$  on a 
curved quantum guide in $\mathbb{R}^n$ ($n=2,3$) with an asymptotically straight reference curve. 
We give  uniqueness results for the inverse problem associated to the reconstruction of  the curvature by using either  observations of spectral data or a  boot-strapping method.\\
 keywords: Inverse Problem, Quantum Guide, Curvature}

\section{Introduction and main results in dimension $n=2$}

\label{sec.1}
The spectral properties of curved quantum guides have been studied intensively for several years, because of their  
applications in quantum mechanics, electron motion. We can cite among several papers \cite{ES}, \cite{GJ}, \cite{KK}, 
\cite{DE}, \cite{DEK}, \cite{CDFK} \ldots \\
However, inverse problems associated with curved quantum guides have not been studied to our knowledge, except in 
\cite{C}. Our aim is to establish  uniqueness results for the inverse problem of  the reconstruction of the curvature of the quantum guide: the data of one eigenpair determines uniquely the curvature up to its sign and  similar results are obtained by considering the knowledge of a solution of Poisson's equation in the guide.\\
We consider the Laplacian operator on a non trivially curved quantum guide $\Omega \subset 
\mathbb{R}^2$ which is not self-intersecting, with Dirichlet boundary conditions, denoted by 
$-\Delta_{D}^{\Omega}.$ We proceed as in \cite{ES}.
We denote by $\Gamma=(\Gamma_1,\Gamma_2)$ the function $C^3$-smooth (see \cite[Remark 5]{CDFK}) which characterizes 
the reference curve and by $N=(N_1,N_2)$ the  outgoing normal to the boundary of $\Omega$.
We denote by $d$ the fixed width of $\Omega$ and by 
$\Omega _0:=\mathbb{R} \times ]-d/2,d/2[.$
Each point $(x,y)$ of $\Omega$ is described by the curvilinear coordinates 
$(s,u)$ as follows:
\begin{equation}\label{f}
\hat{f}: \Omega_0 \longrightarrow \Omega\;\; \mbox{ with }\;\; (x,y)= 
\hat{f}(s,u)=\Gamma(s)+u N(s).
\end{equation}
We assume $\Gamma'_1(s)^2+\Gamma'_2(s)^2=1$ and we recall that the signed 
curvature $\gamma$ of $\Gamma$ is defined by:
\begin{equation}\label{gamma}
\gamma(s)=-\Gamma''_1(s)\Gamma'_2(s)+\Gamma''_2(s)\Gamma'_1(s),
\end{equation}
named so because $|\gamma (s)|$ represents the curvature of the reference 
curve at $s$. We recall that a guide is called simply-bent if $\gamma$ does 
not change sign in $\mathbb{R}$. We assume throughout this article that:\\

\begin{assumption} \label{ass1}
\begin{description}
\item{i)}
$\hat{f}$ is injective.
\item{ii)}
$\gamma \in C^2(\mathbb{R})\cap L^{\infty}(\mathbb{R}),\; \gamma \not \equiv 0,\;
$ (i.e. $\Omega$ is non-trivially curved).
\item{iii)}
$\frac{d}{2}<\frac{1}{ \| \gamma \|_{\infty}},\;  $ where $\| \gamma \|_{\infty}:= 
\sup _{s\in \mathbb{R}} |\gamma (s)| =\| \gamma \|_{L^{\infty}
(\mathbb{R})}.$
\item{iv)}
$\gamma(s) \rightarrow 0 
\mbox{ as } |s|\rightarrow +\infty$ (i.e. $\Omega$ is asymptotically straight).
\end{description}
\end{assumption}

Note that, by the inverse function theorem, the map $\hat{f}$ 
(defined by (\ref{f})) is a local diffeomorphism 
provided $1-u\gamma(s) \neq 0$, for all $u,s,$ which is guaranteed by
Assumption \ref{ass1} and since $\hat{f}$ is assumed to be injective, the map $\hat{f}$ is a global diffeomorphism. Note also that $1-u \gamma(s)>0$ for all $u$ and $s$. (More 
precisely, $0<1-\frac{d}{2} \|\gamma \|_{\infty} \leq 1-u \gamma (s) \leq 
1+\frac{d}{2} \|\gamma \|_{\infty}$ for all $u,s.$)
The curvilinear coordinates $(s,u)$ are locally orthogonal, so by virtue of the Frenet-Serret formulae, the metric in 
$\Omega$ is expressed with respect to them through a 
diagonal metric tensor (e.g. \cite{KK})
\begin{equation}\label{metric}
(g_{ij})= \left ( \begin{array}{cc} (1-u\gamma(s))^2 & 0 \\ 0 & 1 
\end{array} \right ).
\end{equation}
The transition to the curvilinear coordinates represents an isometric map of 
$L^2(\Omega)$ to $L^2(\Omega_0, g^{1/2}\; ds du)$ where
\begin{equation}\label{g}
(g(s,u))^{1/2}: = 1-u \gamma(s)
\end{equation}
is the Jacobian $\displaystyle \frac{\partial(x,y)}{\partial(s,u)}$. 
So we can replace the Laplacian operator $-\Delta_{D}^{\Omega}$ acting on 
$L^2(\Omega)$ by the Laplace-Beltrami operator $H_g$ acting on 
$L^2(\Omega _0,g^{1/2}dsdu)$
relative to the given metric tensor $(g_{ij})$ ( see (\ref{metric}) and 
(\ref{g})) where:
\begin{equation}\label{Hg}
H_g := - g^{-1/2} \partial _s(g^{-1/2}\partial _s) -  
g^{-1/2} \partial _u(g^{1/2}\partial _u).
\end{equation}
We rewrite $H_g$ (defined by (\ref{Hg})) into a Schr\"{o}dinger-type operator 
acting on 
$L^2(\Omega_0, ds du)$.  Indeed, using the unitary transformation
\begin{equation}\label{Ug}
\begin{array}{ccc} U_g : L^2(\Omega_0, g^{1/2}\; ds du) & \longrightarrow  & 
L^2(\Omega_0, ds du) \\
\psi & \mapsto & g^{1/4} \psi \end{array} 
\end{equation}
setting 
$$H_{\gamma}:= U_g H_g U_g^{-1},$$ 
we get 
\begin{equation}\label{Hgamma2}
H_{\gamma}= -\partial_s(c_{\gamma}(s,u) \partial_s)-\partial_u^2+
V_{\gamma}(s,u)
\end{equation}
with
\begin{equation}\label{cgamma}
c_{\gamma}(s,u)=\frac{1}{(1-u \gamma(s))^2}
\end{equation}
and
\begin{equation}\label{Vgamma}
V_{\gamma}(s,u)=-\frac{\gamma^2(s)}{4(1-u \gamma(s))^2}-
\frac{u \gamma''(s)}{2(1-u \gamma(s))^3}-\frac{5u^2 \gamma'^2(s)}
{4(1-u \gamma(s))^4}.
\end{equation}

We will assume throughout all this paper that  the following assumption 
is satisfied:\\

\begin{assumption}\label{ass2} 
$\gamma \in C^2(\mathbb{R}) \mbox{ and } 
\gamma ^{(k)} \in L^{\infty}(\mathbb{R}) \mbox{ for each } k=0,1,2$
 where $\gamma ^{(k)}$ denotes the $k^{th}$ derivative of $\gamma.$
\end{assumption}

{\bf Remarks:}
Since $\Omega$ is non trivially-curved and asymptotically straight, the 
operator $-\Delta_{D}^{\Omega}$ has at least one eigenvalue of finite 
multiplicity below its essential spectrum (see \cite{CDFK}, \cite{KK} ; see also 
\cite{ES} under the 
additional assumptions that the width $d$ is sufficiently small and the 
curvature $\gamma$ is rapidly  decaying at infinity ; see \cite{GJ} under the 
assumption that the curvature $\gamma$ has a compact support).\\
Furthermore, note that such operator $H_{\gamma}$ admits bound states and 
that the minimum  eigenvalue $\lambda_1$ is simple and associated with a 
positive eigenfunction $\phi_1$ (see \cite[Sec.8.17]{GT}). Then, note that by 
\cite[Theorem 7.1]{S} any eigenfunction of $H_{\gamma}$ is continuous  and by 
\cite[Remark 25 p.182]{B} any eigenfunction of $H_{\gamma}$ belongs to 
$H^2(\Omega _0).$ \\
Finally, note also that $(\lambda,\phi)$ is an eigenpair (i.e. an 
eigenfunction associated with its eigenvalue) of the operator $H_{\gamma}$ 
acting on $L^2(\Omega_0,dsdu)$ means that $(\lambda,U_{g}^{-1} \phi)$ is an 
eigenpair of $-\Delta_{D}^{\Omega}$ acting on $L^2(\Omega).$ 
So the data of 
one eigenfunction of the operator $H_{\gamma}$ is equivalent to the data of 
one eigenfunction of $-\Delta_{D}^{\Omega}.$

We first prove that the data of one eigenpair determines 
uniquely the curvature.\\

\begin{theorem}
\label{th1}
Let $\Omega$ be the curved guide in $\mathbb{R}^2$ defined as above. Let 
$\gamma$ be the signed curvature defined by (\ref{gamma}) and satisfying 
Assumptions \ref{ass1}, \ref{ass2}. Let $H_{\gamma}$ be the operator defined by (\ref{Hgamma2}) and 
$(\lambda,\phi)$ be an eigenpair of $H_{\gamma}.$\\
Then
$$ \gamma^2(s) =-4 \frac{\Delta \phi (s,0)}{\phi(s,0)}-4\lambda
$$
 for all 
$s$ when $\phi(s,0)\neq 0$. 
\end{theorem}

Note that the condition $\phi(s,0)\neq 0$ in Theorem \ref{th1} is satisfied for the positive eigenfunction 
$\phi_1$ and for all $s \in \mathbb{R}.$
Then, we prove under \\

\begin{assumption}\label{ass3} 
$\gamma \in C^5(\mathbb{R}) \mbox{ and } \gamma^{(k)} \in 
L^{\infty}(\mathbb{R}) \mbox{ for each } k=0,\dots,5,$
\end{assumption}
  
\noindent that one weak solution $\phi$ of the problem 
\begin{equation}\label{P}
\left \{
\begin{array}{ll}
H_{\gamma}\phi =f \mbox{ in } \Omega_0 \\
\phi=0 \mbox{ on } \partial \Omega_0
\end{array} \right.
\end{equation}
 (where $f$ is a known given function) is in fact a classical solution and the data of $\phi$ determines uniquely 
the curvature $\gamma$.

\begin{theorem}
\label{th2}
Let $\Omega$ be the curved guide in $\mathbb{R}^2$ defined as above. Let 
$\gamma$ be the signed curvature defined by (\ref{gamma}) and satisfying 
Assumptions \ref{ass1} and \ref{ass3}. Let $H_{\gamma}$ be the operator defined by (\ref{Hgamma2}). Let $f 
\in H^3(\Omega_0) \cap C(\Omega_0)$ and let 
$\phi \in H^1_{0}(\Omega_0)$ be a weak solution of (\ref{P}).\\
Then we have $\displaystyle \gamma^2(s) =-4 \frac{\Delta \phi (s,0)}{\phi(s,0)}-
4\frac{f(s,0)}{\phi (s,0)}$ for all $s$ when $\phi(s,0)\neq 0 $
\end{theorem}

In the case of a simply-bent guide (i.e. when $\gamma$ does not change sign 
in $\mathbb{R}$), we can restrain the hypotheses upon the regularity of 
$\gamma.$ 
We obtain the following result:

\begin{theorem}
\label{th3}
Let $\Omega$ be the curved guide in $\mathbb{R}^2$ defined as above. Let 
$\gamma$ be the signed curvature defined by (\ref{gamma}) and satisfying 
Assumptions \ref{ass1} and \ref{ass2}. 
We assume also that $\gamma$ is a nonnegative function. Let 
$H_{\gamma}$ be the operator defined by (\ref{Hgamma2}). Let  
$f \in L^2(\Omega_0)$ 
be a non null function and let $\phi$ be a weak solution in $H^1_0(\Omega_0)$ of (\ref{P})
Assume that there exists a positive constant $M$ such that 
$|f(s,u)|\leq M| \phi (s,u)|$ almost everywhere in 
$\Omega_0.$ Then $(f,\phi)$ 
determines uniquely the curvature $\gamma.$
\end{theorem}

Note that the above result is still valid for a nonpositive function $\gamma$.\\

This paper is organized as follows: In Section 2, we prove Theorems \ref{th1}, \ref{th2} and \ref{th3}. In 
Sections 3 and 4, we extend our results to the case of a curved quantum guide defined in $\mathbb{R}^3.$

\section{Proofs of Theorems \ref{th1}, \ref{th2} and \ref{th3}}
\label{sec.2}
\subsection{Proof of Theorem \ref{th1}}
Recall that $\phi$ is an eigenfunction of $H_{\gamma},$ belonging to
$H^2(\Omega_0).$  Since $\phi$ is 
continuous and $H_{\gamma} \phi =\lambda \phi,$ then $H_{\gamma} \phi$ is
continuous too. Thus, noticing that $c_{\gamma}(s,0)=1,$ we
deduce the continuity of the function $(s,0) \mapsto \Delta \phi (s,0)$ and
from (\ref{Hgamma2}) to (\ref{Vgamma}), we get:
$$-\Delta \phi(s,0)- \frac{\gamma^{2}(s)}{4} 
\phi (s,0) =\lambda \phi(s,0)
$$
 and equivalently, 
$$\gamma ^{2}(s) =-4 \frac{
\Delta \phi (s,0)}{\phi(s,0)}-4\lambda  \mbox{ if } \phi(s,0)\neq 0.$$

\subsection{Proof of Theorem \ref{th2}}
First, we recall from \cite[Remark 25 p.182]{B} the following lemma.

\begin{lemma}\label{Lem}
For a second-order elliptic operator defined in a domain $\omega \subset 
\mathbb{R}^n$, if $\phi \in H^1_0(\omega)$ satisfies 
$$\int_{\omega} \sum_{i,j} a_{ij} \frac{\partial \phi}{\partial x_i}
\frac{\partial \psi}{\partial x_j}=\int_{\omega}f \psi \mbox{ for all }
\psi \in H^1_0(\omega)$$ 
then if $\omega$ is of class $C^2$
$$(f \in L^2(\omega), \, a_{ij} \in C^1(\overline{\omega}),\, D^{\alpha} 
a_{ij} 
\in L^{\infty}(\omega) \mbox{ for all } i,j \mbox{ and for all } \alpha, 
|\alpha| \leq 1) $$
$$\mbox{ imply } (\phi \in H^2(\omega))$$
and for $m \geq 1,$ if $\omega$ is of class $C^{m+2}$
$$(f \in H^m(\omega), \, a_{ij} \in C^{m+1}(\overline{\omega}),\, D^{\alpha} 
a_{ij} \in 
L^{\infty}(\omega) \mbox{ for all } i,j \mbox{ and for all } \alpha, 
|\alpha| \leq m+1) $$
$$\mbox{ imply } (\phi \in H^{m+2}(\omega)).$$
\end{lemma}

Now we can prove the Theorem \ref{th2}.\\

We have $H_{\gamma} \phi=f$, so  
\begin{equation}
\label{21}
\int_{\Omega_0}[c_{\gamma}(\partial_{s} \phi)(\partial_{s} \psi)+(\partial_{u}
 \phi)(\partial_{u} \psi)]=\int_{\Omega_0}[f -V_{\gamma}\phi]\psi 
\mbox{ for all } \psi \in H^1_0 (\Omega_0)
\end{equation}
with $c_{\gamma}$ defined by (\ref{cgamma}) and $V_{\gamma}$ defined by 
(\ref{Vgamma}).\\

Using Assumption \ref{ass3}, since $\gamma^{(k)} \in L^{\infty}(\Omega_0)$ 
for $k=0,1,2$ then $V_{\gamma} \in L^{\infty}(\Omega_0)$ and $f -
V_{\gamma} \phi \in L^2(\Omega_0).$ From the hypotheses $\gamma \in 
C^1(\mathbb{R})$ and $\gamma' \in L^{\infty}(\mathbb{R})$, we get that 
$c_{\gamma} \in C^1(\overline{\Omega_0}),\, D^{\alpha} c_{\gamma} \in 
L^{\infty}
(\Omega_0)$ for any $\alpha, \, |\alpha |\leq 1,$ and so, using Lemma \ref{Lem} 
for the equation (\ref{21}), we 
obtain that $\phi \in H^2(\Omega_0).$\\
By the same way, we get that $f -V_{\gamma}\phi \in H^1(\Omega_0),
\, c_{\gamma} \in C^2(\overline{\Omega_0})$ and $D^{\alpha}c_{\gamma} \in 
L^{\infty}(\Omega_0)$ for any $\alpha,\, |\alpha |\leq 2$ (from $\gamma 
\in C^3(\mathbb{R}),\, \gamma ^{(k)} \in L^{\infty}(\mathbb{R})$ for any 
$k=0,\dots,3$). Using Lemma \ref{Lem}, we obtain that $\phi \in H^3(\Omega_0).$\\
We apply again the Lemma \ref{Lem} to get that $\phi \in H^4(\Omega_0)$ (since 
$f -V_{\gamma}\phi \in H^2(\Omega_0),\, c_{\gamma}\in 
C^3(\overline{\Omega_0}),\, D^{\alpha}c_{\gamma} \in L^{\infty}(\Omega_0)$ for all 
$\alpha,\, |\alpha|\leq 3$, from the hypotheses 
$\gamma \in C^4(\mathbb{R})$ and 
$\gamma^{(k)} \in L^{\infty}(\mathbb{R})$ for $k=0,\dots,4.$).\\
Finally, using Assumption \ref{ass3} and Lemma \ref{Lem}, we obtain that $\phi \in 
H^5(\Omega_0)$.

Due to the regularity of $\Omega_0,$ we have $\phi \in H^5(\mathbb{R}^2)$
and $\Delta \phi  \in H^3(\mathbb{R}^2).$ 
Since $\nabla (\Delta \phi) \in (H^2(\mathbb{R}^2))^2$ and 
$H^2(\mathbb{R}^2)\subset L^{\infty}(\mathbb{R}^2)$, we can deduce that $\Delta \phi$ is continuous 
(see \cite[Remark 8 p.154]{B}).\\ 
Therefore we can conclude by using the continuity of the function 
$$(s,0) \mapsto -\partial_s
(c_{\gamma}(s,0) \partial_s \phi (s,0)) -\partial_u ^2
\phi(s,0)=f(s,0)-V_{\gamma} (s,0)\ \phi(s,0).$$
Therefore, we get: $-\Delta \phi(s,0)- \frac{\gamma^{2}(s)}{4} 
\phi (s,0) =f(s,0)$ and equivalently, 
$$\gamma ^{2}(s) =-4 \frac{
\Delta \phi (s,0)}{\phi(s,0)}-4\frac{ f(s,0)}{\phi (s,0)} \mbox{ if } 
\phi(s,0)\neq 0.$$

\subsection{Proof of Theorem \ref{th3}}

We prove here that $(f,\phi)$ determines uniquely $\gamma$ when 
$\gamma$ is a nonnegative function.\\
For that, assume that $\Omega _1$ and $\Omega _2$ are two quantum guides 
in $\mathbb{R}^2$ with same width $d$. We denote by $\gamma_1$ and $\gamma_2$ 
the curvatures respectively associated with $\Omega_1$ and $\Omega_2$ and we 
suppose that each $\gamma_i$ satisfies Assumption \ref{ass2} and is a nonnegative 
function.
Assume that $H_{\gamma_1} \phi = f = H_{\gamma_2} \phi$.\\
Then $\phi$ satisfies 
\begin{equation} \label{(2.2)}
 -\partial_s((c_{\gamma_1}(s,u)-c_{\gamma_2}(s,u)) \partial_s \phi(s,u))+
(V_{\gamma_1}(s,u)-V_{\gamma_2}(s,u))\phi(s,u)=0.
\end{equation}
Assume that $\gamma_1 \not\equiv \gamma_2.$\\
Step 1. First, we consider the case where (for example) 
$\gamma_1 (s) < \gamma_2 (s)$ for all $s \in \mathbb{R}.$\\
Let $\epsilon >0,\; \omega_{\epsilon}:=\mathbb{R} \times I_{\epsilon}$ with 
$I_{\epsilon}=]-\epsilon,0[.$ 
Multiplying (\ref{(2.2)}) by $\phi$ and integrating over ${\omega_{\epsilon}}$, 
we get: 
\begin{equation}\label{(2.3)}
\int_{\omega_{\epsilon}} (c_{\gamma_1}-c_{\gamma_2}) (\partial_s \phi)^2 -
\int_{\partial \omega_{\epsilon}}(c_{\gamma_1}-c_{\gamma_2}) 
(\partial_s \phi) \phi \nu_s +
\int_{\omega_{\epsilon}}(V_{\gamma_1}-V_{\gamma_2})\phi^2=0.
\end{equation}
Since $\epsilon <<1$, $ V_{\gamma_i}(s,u) \simeq -\frac{\gamma_i^2(s)}{4}$ 
for $i=1,2$, and so 
$V_{\gamma_1}(s,u) -V_{\gamma_2}(s,u)>0$ in $\omega_{\epsilon}$.\\
Moreover, since
\begin{equation}\label{(2.4)}
c_{\gamma_1}(s,u)-c_{\gamma_2}(s,u)=\frac{u(\gamma_1(s)-\gamma_2(s))
(2-u(\gamma_1(s)+\gamma_2(s))}
{(1-u\gamma_1(s))^2(1-u\gamma_2(s))^2},
\end{equation}
we have $c_{\gamma_1}(s,u)> c_{\gamma_2}(s,u)$ in $\omega_{\epsilon}$.\\
Since
\begin{equation}\label{(2.5)}
\int_{\partial \omega_{\epsilon}}(c_{\gamma_1}-c_{\gamma_2}) 
(\partial_s \phi) \phi \nu_s =0,
\end{equation}
Thus from (\ref{(2.3)})-(\ref{(2.5)}), we get
\begin{equation}\label{(2.6)}
\int_{\omega_{\epsilon}} (c_{\gamma_1}-c_{\gamma_2}) (\partial_s \phi)^2 +
\int_{\omega_{\epsilon}}(V_{\gamma_1}-V_{\gamma_2})\phi^2= 0
\end{equation}
with $c_{\gamma_1}-c_{\gamma_2} > 0 \mbox{ in } \omega_{\epsilon}$ and 
$V_{\gamma_1}-V_{\gamma_2}>0
\mbox{ in } \omega_{\epsilon}.$
We can deduce that $\phi=0$ in $\omega_{\epsilon}$.\\
Using a unique continuation theorem (see \cite[Theorem XIII.63 p.240]{RS}), 
from 
$H_{\gamma}\phi =f$, noting that $-\Delta(U_{g}^{-1}\phi) =
 U_{g}^{-1} f= g^{-1/4}f,$ 
(recall that $U_g$ is defined by (\ref{Ug})) 
and so by $|f| \leq M | \phi |$ we have 
$\vert \Delta (U_g ^{-1}\phi)\vert \leq M \vert g^{-1/4} \phi \vert \
 \mbox{ with }  g>0 \, \mbox{ a.e., }$ and 
we can deduce that $\phi=0$ in $\Omega_0$. 
So we get a contradiction (since $H_{\gamma }\phi=f$ and $f$ is assumed to be 
a non null function).\\

Step 2. From Step 1, we obtain that there  exists at least one point $s_0 \in 
\mathbb{R}$ such that $\gamma_1(s_0)=\gamma_2(s_0).$ Since $\gamma_1 \not\equiv 
\gamma_2$, we can choose 
 $a\in \mathbb{R}$ and $b\in \mathbb{R} \cup \{+\infty\}$ such that
(for example) $\gamma_1(a) = \gamma_2(a),\; \gamma_1(s) < \gamma_2(s)
 \mbox{ for all } s \in ]a,b[ \mbox{ and } 
\gamma_1(b) = \gamma_2(b) \mbox{ if } b \in \mathbb{R}.$\\
We proceed as in Step 1, considering, in this case, 
$\omega_{\epsilon}:= ]a,b[ \times I_{\epsilon}$. We study again the 
equation (\ref{(2.3)}) and as in Step 1, we have
$$\int_{\partial \omega_{\epsilon}}(c_{\gamma_1}-c_{\gamma_2}) 
(\partial_s \phi) \phi \nu_s =0.$$
Indeed from (\ref{(2.4)}) and $\gamma_1 (a)=\gamma_2 (a)$ we have 
$c_{\gamma_1} (a,u)=
c_{\gamma_2} (a,u)$ and so 
$$\int_{-\epsilon} ^0 (c_{\gamma_1} (a,u)-c_{\gamma _2} (a,u)) \partial_s 
\phi(a,u) \phi(a,u) \ du =0.$$
 By the same way
if $b\in \mathbb{R}$, we also have  
$c_{\gamma_1} (b,u)=c_{\gamma_2} (b,u).$ 
Thus the equation (\ref{(2.3)}) becomes (\ref{(2.6)})
with $c_{\gamma_1}-c_{\gamma_2} >0 \mbox{ in } \omega_{\epsilon}$ and 
$V_{\gamma_1}-V_{\gamma_2}>0 \mbox{ in } \omega_{\epsilon}.$
So $\phi=0$ in $\omega_{\epsilon}$ and as in Step 1, by a unique 
continuation theorem, we obtain that $\phi=0$ in $\Omega_0$. 
Therefore we get a contradiction.\\

Note that the previous theorem is true if we replace the 
hypothesis 
$"\gamma$ is nonnegative" by the hypothesis $"\gamma$ is nonpositive". 
Indeed, in this last case, 
we just have to take $I_{\epsilon}=]0,\epsilon[$ 
and the proof rests valid.

\section{Uniqueness result for a $\mathbb{R}^3$-quantum guide}
Now, we apply the same ideas for a tube $\Omega$ in $\mathbb{R}^3.$
We proceed here as in \cite{CDFK}. 
Let $s\mapsto \Gamma (s),\, \Gamma=(\Gamma_1, 
\Gamma_2, \Gamma_3),$ be a curve in $\mathbb{R}^3.$
We assume that $\Gamma : \mathbb{R} \rightarrow \mathbb{R}^3$ is a $C^4$-smooth curve satisfying the following hypotheses
\begin{assumption}\label{A3}
$\Gamma$ possesses a positively oriented Frenet frame $\{e_1,e_2,e_3\}$ with the properties that
\begin{description}
\item{i)}
$e_1=\Gamma'$,
\item{ii)}
$\forall i  \in \{1,2,3\}, \; e_i \in C^1(\mathbb{R}, \mathbb{R}^3),$
\item{iii)}
$\forall i  \in \{1,2\}, \; \forall s \in \mathbb{R},\; e_i'(s)$ lies in the span of $e_1(s), \ldots, e_{i+1}(s)$.
\end{description}
\end{assumption}
Recall that a sufficient condition to ensure the existence of the Frenet frame of Assumption \ref{A3} is to require 
that for all $s \in \mathbb{R}$ the vectors $\Gamma'(s),\Gamma''(s)$ are linearly independent.\\
Then we define the moving frame $\{\tilde{e}_1, \tilde{e}_2, \tilde{e}_3\}$ along $\Gamma$ by following \cite{CDFK}. This moving frame better reflects the geometry of the curve and it is still called the Tang frame because it is a generalization of the Tang frame known from the theory of three-dimensional waveguides.\\
Given a $C^5$  bounded open connected neighborhood 
$\omega$ of $(0,0) \in 
\mathbb{R}^2$, let $\Omega_0$ denote the straight tube $\mathbb{R} \times 
\omega.$ We define the curved tube $\Omega$ of cross-section $\omega$ about 
$\Gamma$ by
\begin{equation}\label{ftilde}
\Omega := \tilde{f}(\mathbb{R} \times \omega)=\tilde{f}(\Omega_0),\, 
\tilde{f}(s,u_2,u_3):= \Gamma (s)+\sum_{i=2}^{3}u_i \sum_{j=2}^{3} R_{ij}(s) e_j(s) = \Gamma (s)+\sum_{i=2}^{3}u_i\tilde{e}_i(s)
\end{equation}
with $u=(u_2,u_3) \in \omega$ and 
$$R(s):=(R_{ij}(s))_{i,j\in \{2,3\}} =\left ( \begin{array}{cc} 
\cos (\theta (s)) & - \sin (\theta (s)) \\
\sin (\theta (s)) & \cos (\theta (s))
\end{array} \right ),$$
$\theta$ being a real-valued differentiable function such that 
$\theta'(s) =\tau (s)$ the torsion of $\Gamma.$ This differential equation is a consequence of the definition of the moving Tang frame (see \cite[Remark 3] {CDFK}).\\
Note that $R$ is a rotation matrix in $\mathbb{R}^2$ chosen in such a way 
that $(s,u_2,u_3)$ are orthogonal ``coordinates'' in $\Omega$.
Let $k$ be the first curvature function of $\Omega.$ Recall that since 
$\Omega \subset \mathbb{R}^3,$ $k$ is a nonnegative function. We assume 
throughout all this section that the following hypothesis holds:\\

\begin{assumption}\label{Ass4} 
\begin{description}
\item {i)} $k \in C^2 (\mathbb{R}) \cap L^{\infty}(\mathbb{R}),\,
a:= \sup_{u\in \omega} \|u\|_{\mathbb{R}^2} <\frac{1}{\|k\|_{\infty}},\, 
k(s) \rightarrow 0 \mbox{ as } |s| \rightarrow +\infty$ 
\item {ii)} $\Omega$ does not overlap.
\end{description}
\end{assumption}

The Assumption \ref{Ass4} assures that the map $\tilde{f}$ (defined by (\ref{ftilde}))
is a diffeomorphism (see \cite{CDFK}) in 
order to identify $\Omega$ with the Riemannian manifold $(\Omega_0, (g_{ij}))$
where $(g_{ij})$ is the metric tensor induced by $\tilde{f}$, i.e. 
$(g_{ij}):= ^{t} J(\tilde{f}) . J(\tilde{f}),$
($J(\tilde{f})$ denoting the Jacobian matrix of $\tilde{f}$). Recall that $(g_{ij})=
diag (h^2,1,1)$ (see [3]) with 
\begin{equation} \label{h}
h(s,u_2,u_3):= 1- k(s)(\cos (\theta (s)) u_2 +\sin (\theta (s)) u_3).
\end{equation}
Note that Assumption \ref{Ass4} implies that $0<1-a\|k\|_{\infty} \leq 1-
h(s,u_2,u_3)\leq 1+a \|k\|_{\infty}$ for all $s \in \mathbb{R}$ and 
$u=(u_2,u_3)\in \omega.$ Moreover, setting 
\begin{equation}
g:=h^2
\end{equation}
we can replace the Dirichlet Laplacian operator 
$-\Delta_{D}^{\Omega}$ acting on $L^2(\Omega)$ by the Laplace-Beltrami 
operator $K_g$ acting on $L^2(\Omega_0, hdsdu)$ relative to the metric 
tensor $(g_{ij}).$ We can rewrite $K_g$ into a Schr\"{o}dinger-type operator 
acting on $L^2(\Omega_0, dsdu)$.
Indeed, using the unitary transformation
\begin{equation}\label{Wg}
\begin{array}{ccc} 
W_g : L^2(\Omega_0, hdsdu) & \longrightarrow  & L^2(\Omega_0, dsdu) \\
\psi & \mapsto & g^{1/4} \psi \end{array} 
\end{equation}
setting 
\begin{equation}\label{Hk1}
H_{k}: = W_g K_g W_g^{-1},
\end{equation}
we get 
\begin{equation}\label{Hk2}
H_{k}= -\partial_s(h^{-2}\partial_s)
-\partial_{u_2}^2 -\partial_{u_3}^2+V_{k}
\end{equation}
where $\partial_s$ denotes the derivative relative to $s$ and $\partial_{u_i}$
denotes the derivative relative to $u_i$ and with
\begin{equation}\label{Vk}
V_{k}:=-\frac{k^2}{4h^2}+\frac{\partial_s^2 h}{2h^3}
-\frac{5(\partial _s h)^2}{4h^4}.
\end{equation}
We assume also throughout all this section that the following hypotheses hold:\\

\begin{assumption}\label{Ass5} 
\begin{description}
\item {i)} $k'\in L^{\infty}(\mathbb{R}),\, k'' \in 
L^{\infty}(\mathbb{R})$
\item {ii)} $\theta \in C^2(\mathbb{R}),\, 
\theta'=\tau \in L^{\infty}(\mathbb{R}),\,
\theta'' \in L^{\infty}(\mathbb{R}).$
\end{description}
\end{assumption}

{\bf Remarks:} Note that, as for the 2-dimensional case, such operator 
$H_{k}$ (defined by (\ref{h})-(\ref{Vk})) admits bound states and that the minimum 
eigenvalue $\lambda_1$ is simple and associated with a positive eigenfunction
 $\phi_1$ (see \cite{CDFK,GT}). Still note that $(\lambda,\phi)$ is an eigenpair 
of the operator $H_k$ acting on $L^2(\Omega_0, dsdu)$ means that 
$(\lambda, W_g ^{-1} \phi)$ is an eigenpair of $-\Delta _{D}^{\Omega}$ acting 
on $L^2(\Omega)$ (with $W_g$ defined by (\ref{Wg})). Finally, note that by \cite[Theorem 7.1]{S} 
any eigenfunction of $H_k$ is continuous and by \cite[Remark 25 p.182]{B} 
any eigenfunction of $H_k$ belongs to $H^2(\Omega_0).$\\

As for the 2-dimensional case, first we prove that the data of one eigenpair 
determines uniquely the curvature.\\

\begin{theorem}\label{th4}
Let $\Omega$ be the curved guide in $\mathbb{R}^3$ defined as above. Let 
$k$ be the first curvature function of $\Omega$. Assume that Assumptions \ref{A3} to \ref{Ass5} are satisfied. 
Let $H_k$ be the operator defined by (\ref{h})-(\ref{Vk}) and $(\lambda,\phi)$ be an eigenpair of $H_k$. \\
Then $k^2(s) =-4 \frac{\Delta \phi (s,0,0)}{\phi(s,0,0)}-4\lambda$ for all 
$s$ when $\phi(s,0,0)\neq 0$. 
\end{theorem}

Then, under 

\begin{assumption}\label{Ass6} 
\begin{description}
\item {i)} $k \in C^5(\mathbb{R}), k^{(i)} \in L^{\infty}(\mathbb{R}) 
\mbox{ for all } i=0,\dots,5$
\item {ii)} $\theta \in C^5(\mathbb{R}), \theta^{(i)} \in 
L^{\infty}(\mathbb{R}) \mbox{ for all } i=1,\dots,5$
\end{description}
\end{assumption}

where $k^{(i)}$ (resp. $\theta^{(i)}$) denotes the i-th derivative of $k$ 
(resp. of $\theta$), we obtain the following result:

\begin{theorem}\label{th5}
Let $\Omega$ be the curved guide in $\mathbb{R}^3$ defined as above. Let 
$k$ be the first curvature function of $\Omega$. Assume that Assumptions \ref{A3} to \ref{Ass6} are satisfied. 
Let $H_k$ be the operator defined by (\ref{h})-(\ref{Vk}). Let $f \in H^3(\Omega _0) \cap 
C(\Omega_0)$ and let $\phi \in H^1_0(\Omega _0)$ be a weak solution of $H_k \phi=f$ in $\Omega_0.$ \\
Then  $\phi$ is a classical solution and $k^2(s) =-4 \frac{\Delta \phi (s,0,0)}{\phi(s,0,0)}-4\frac{ f(s,0,0)}
{\phi(s,0,0)}$ for all $s$ when $\phi(s,0,0)\neq 0$. 
\end{theorem}

{\bf Remarks:} Recall that in $\mathbb{R}^3$, $k$ is a nonnegative function 
and that the condition imposed on $\phi$ ($\phi (s,0,0) \neq 0$) in 
Theorems \ref{th4} and \ref{th5} is satisfied by the positive eigenfunction $\phi_1.$\\

As for the two-dimensional case, we can restrain the hypotheses upon the 
regularity of the functions $k$ and $\theta$. \\

For a guide with a known torsion, we obtain the following result:

\begin{theorem}\label{th6}
Let $\Omega$ be the curved guide in $\mathbb{R}^3$ defined as above. Let 
$k$ be the first curvature function of $\Omega$ and let $\tau$ be the second 
curvature function (i.e. the torsion) of $\Omega$. Denote by $\theta$ a primitive  
of $\tau$ and suppose that $0\leq \theta(s) \leq \frac{\pi}{2}$ for all $s\in \mathbb{R}.$
Assume that Assumptions \ref{A3} to \ref{Ass5} are satisfied. Let $H_k$ 
be the operator defined by (\ref{h})-(\ref{Vk}). Let $f \in L^2(\Omega_0)$ be a non 
null function and let $\phi \in H^1_0(\Omega_0)$ be a weak solution of $H_k \phi =f$ in $\Omega_0.$ Assume that 
there exists a positive constant $M$ such that $|f(s,u)| \leq M | \phi (s,u)|$ almost everywhere in $\Omega_0.$\\
Then the data $(f,\phi)$ determines uniquely the first curvature function $k$ if the torsion $\tau$ is given.
\end{theorem}

\section{Proofs of Theorem \ref{th4}, \ref{th5} and \ref{th6} }
\label{sec.3}
\subsection{Proof of Theorem \ref{th4}}

Recall that $\phi$ is an eigenfunction of $H_k.$ Since $\phi$ is continuous, 
$H_k \phi =\lambda \phi$ and $\phi \in H^2(\Omega_0)$ then $H_k \phi$ is 
continous. 
Therefore, for $u=(u_2,u_3)=(0,0)$, we get: 
$-\Delta \phi(s,0,0)- \frac{k^{2}(s)}{4} 
\phi (s,0,0) =\lambda \phi(s,0,0)$ and equivalently, 
$k^{2}(s) =-4 \frac{\Delta \phi (s,0,0)}{\phi(s,0,0)}-4\lambda$ if 
$\phi(s,0,0)\neq 0.$ 

\subsection{Proof of Theorem \ref{th5}}

We follow the proof of Theorem \ref{th2}. 
We have $H_{k} \phi=f$ with $\phi \in H^1_0(\Omega_0)$. So
\begin{equation}\label{4.1}
\int_{\Omega_0}[h^{-2}(\partial_{s} \phi)(\partial_{s} \psi)+(\partial_{u_2}
 \phi)(\partial_{u_2} \psi)+(\partial_{u_3}\phi)(\partial_{u_3} \psi)
]=\int_{\Omega_0}[f -V_{k}\phi]\psi 
\mbox{ for all } \psi \in H^1_0 (\Omega_0)
\end{equation}
with $h$ defined by (\ref{h}) and $V_{k}$ defined by (\ref{Vk}).\\

From Assumptions \ref{Ass4} and \ref{Ass5}, since $k,k',k'',\theta',\theta''$ are 
bounded, we deduce that $V_{k} \in L^{\infty}(\Omega_0)$. Therefore 
$f -V_{k} \phi \in L^2(\Omega_0).$ Moreover  we have also $h^{-2}
\in C^1(\overline{\Omega_0})$ and $D^{\alpha} (h^{-2}) \in L^{\infty}
(\Omega_0)$ for any $\alpha, \, |\alpha |\leq 1.$ Thus, using Lemma \ref{Lem} 
for the equation (\ref{4.1}), we 
obtain that $\phi \in H^2(\Omega_0).$\\
By the same way, we get that $f -V_{k}\phi \in H^1(\Omega_0),
\, h^{-2} \in C^2(\overline{\Omega_0})$ and $D^{\alpha}(h^{-2}) \in 
L^{\infty}(\Omega_0)$ for any $\alpha,\, |\alpha |\leq 2$ (since 
$k \in C^3(\mathbb{R}),\, \theta \in C^3(\mathbb{R})$ and all of their
derivatives are bounded). Using Lemma \ref{Lem}, we obtain that 
$\phi \in H^3(\Omega_0).$\\
We apply again the Lemma \ref{Lem} to get that $\phi \in H^4(\Omega_0)$ (since 
$f -V_{\gamma}\phi \in H^2(\Omega_0),\, c_{\gamma}\in 
C^3(\overline{\Omega_0}),\, D^{\alpha}c_{\gamma} \in L^{\infty}(\Omega_0)$ for all 
$\alpha,\, |\alpha|\leq 3$, from the hypotheses 
$\gamma \in C^4(\mathbb{R})$ and 
$\gamma^{(k)} \in L^{\infty}(\mathbb{R})$ for $k=0,\dots,4.$).\\
Finally, using Assumption \ref{Ass6} and Lemma \ref{Lem}, we obtain that $\phi \in 
H^5(\Omega_0)$. Due to the regularity of $\Omega_0$ (see \cite[Note p.169]{B}), 
we have 
$\phi \in H^5(\mathbb{R}^3)$ and $\Delta \phi  \in H^3(\mathbb{R}^3).$ 
Since $\nabla (\Delta \phi) \in (H^2(\mathbb{R}^3))^3$ and 
$H^2(\mathbb{R}^3)\subset L^{\infty}(\mathbb{R}^3)$, we can deduce that $\Delta \phi$ is continuous 
(see \cite[Remark 8 p.154]{B}).\\ 
Thus we conclude as in Theorem \ref{th2} and for $u=(u_2,u_3)=(0,0)$, we get: 
$-\Delta \phi(s,0,0)- \frac{k^{2}(s)}{4} 
\phi (s,0,0) =f(s,0,0)$ and equivalently, 
$k^{2}(s) =-4 \frac{\Delta \phi (s,0,0)}{\phi(s,0,0)}-4\frac{f(s,0,0)}{
\phi(s,0,0)}$ if $\phi(s,0,0)\neq 0.$

\subsection{Proof of Theorem \ref{th6}}

We prove here that $(f,\phi, \theta)$ determines uniquely $k.$\\
Assume that $\Omega _1$ and $\Omega _2$ are two guides in $\mathbb{R}^3.$ 
We denote by $k_1$ and $k_2$ the first curvatures functions associated with 
$\Omega_1$ and $\Omega_2$ and we denote by $\theta$ a primitive of $\tau$ the 
common torsion of $\Omega_1$ and $\Omega_2.$ We suppose that $k_1, k_2$ and 
$\theta$ satisfy the Assumptions \ref{Ass4} and \ref{Ass5} and that $0\leq \theta (s) \leq \frac{\pi}
{2}$ for all $s \in \mathbb{R}.$ Assume that $H_{k_1} \phi= f =
H_{k_2}\phi.$\\
Then $\phi$ satisfies 
$$-\partial_s ((h_1^{-2}(s,u_2,u_3) -h_2^{-2}(s,u_2,u_3)) \partial_s \phi 
(s,u_2,u_3) )$$
\begin{equation}\label{4.2}
+(V_{k_1}(s,u_2,u_3)-V_{k_2}(s,u_2,u_3)) \phi(s,u_2,u_3)=0
\end{equation}
where $h_1$ (associated with $k_1$) is defined by (\ref{h}), $V_{k_1}$ is 
defined by (\ref{Vk}), $h_2$ (associated with $k_2$) is defined by (\ref{h}) and 
$V_{k_2}$ is defined by (\ref{Vk}).\\
Assume that $k_1 \not\equiv k_2.$\\

Step 1. First, we consider the case where (for example) $k_1 (s) < k_2 (s)$ 
for all $s \in \mathbb{R}.$ Recall that each $k_i$ is a nonnegative 
function.\\
Let $\epsilon>0$ and denote by $J_{\epsilon}:=]-\epsilon,0[ \times 
]-\epsilon,0[,\, O_{\epsilon}:= \mathbb{R} \times J_{\epsilon}$ with $\epsilon$ small enough to have 
$J_{\epsilon}\subset \omega$ (recall that $\Omega_0=\mathbb{R}
\times \omega$).\\
Multiplying (\ref{4.2}) by $\phi$ and integrating over $O_{\epsilon}$, we get:
\begin{equation}\label{4.3}
 \int_{O_{\epsilon}} (h_1^{-2}-h_2^{-2})(\partial_s \phi)^2
+\int_{\partial O_{\epsilon}}(h_1^{-2}-h_2^{-2}) 
(\partial_s \phi) \phi \nu_s
 +\int_{O_{\epsilon}}(V_{k_1}-V_{k_2})\phi^2=0.
\end{equation}
Since $\epsilon<<1,\, V_{k_i}\simeq -\frac{k_i^2(s)}{4}$ 
for $i=1,2$, and so 
$V_{k_1}(s,u_2,u_3) -V_{k_2}(s,u_2,u_3)>0$ in $O_{\epsilon}$.\\
Moreover, note that:
\begin{equation}\label{4.4}
h_1^{-2}(s,u_2,u_3)-h_2^{-2}(s,u_2,u_3)=\frac{\alpha (s,u_2,u_3) 
(k_1(s)-k_2(s))(h_1(s,u_2,u_3)+h_2(s,u_2,u_3))}
{h_1^2 (s,u_2,u_3) h_2^2 (s,u_2,u_3)}
\end{equation}
with $\alpha (s,u_2,u_3):= \cos (\theta (s))u_2+\sin (\theta(s)) u_3.$\\
Since $(u_2,u_3) \in J_{\epsilon}$ and $0\leq \theta (s) \leq \frac{\pi}{2}$ 
for all $s \in \mathbb{R},$ we have $\alpha (s,u_2,u_3)<0.$ Therefore, by 
(\ref{4.4}), we deduce that $h_1^{-2}-h_2^{-2}>0$ in $O_{\epsilon}.$\\
Thus $\int_{O_{\epsilon}} (h_1^{-2}-h_2^{-2})(\partial_s \phi)^2
+\int_{O_{\epsilon}}(V_{k_1}-V_{k_2})\phi^2\geq 0.$\\
Note also that:
\begin{equation}\label{4.5}
\int_{\partial O_{\epsilon}} (h^{-2}_{1}-h^{-2}_{2})
(\partial_s \phi) \phi \nu_s =0.
\end{equation}
Therefore, from (\ref{4.3}) and (\ref{4.5}) we get: 
\begin{equation}\label{4.6}
\int_{O_{\epsilon}} (h^{-2}_{1}-h^{-2}_{2}) (\partial_s \phi)^2
+\int_{O_{\epsilon}}(V_{k_1}-V_{k_2})\phi^2=0
\end{equation}
with $h^{-2}_{1}-h^{-2}_{2}>0 \mbox{ in } O_{\epsilon} \mbox{ and } 
V_{k_1}-V_{k_2}>0 \mbox{ in } O_{\epsilon}.$ \\
From (\ref{4.6}) we can deduce that $\phi=0$ in $O_{\epsilon}.$ Using a unique 
continuation 
theorem (see \cite[Theorem XIII.63 p.240]{RS}), from $H_{k_1}\phi =f,$ 
noting that $-\Delta (W_{g}^{-1} \phi) =W_{g}^{-1} f =
g^{-1/4}f$, by $|f| \leq M |\phi|$ a.e. in  $\Omega_0,$ we can deduce that $\phi=0$ in 
$\Omega_0$. So we get a contradiction since $f$ is assumed to be a non null 
function.\\

Step2. From Step 1, we obtain that there exists at least one point $s_0 \in 
\mathbb{R}$ such that $k_1(s_0)=k_2(s_0).$ Since $k_1 \not\equiv k_2,$ we can 
choose 
$a \in \mathbb{R}$ and $b\in \mathbb{R} \cup \{+\infty\}$ such that (for 
example) $k_1(a)=k_2(a),\, k_1(s)<k_2(s) \mbox{ for all } s\in ]a,b[ 
\mbox{ and } k_1(b)=k_2(b) \mbox{ if } b\in \mathbb{R}.$ We proceed as in 
Step 1, considering in this case $O_{\epsilon}:=]a,b[\times J_{\epsilon}.$ 
From $k_1(a)=k_2(a),$ we get that $h_1^{-2}(a,u_2,u_3)=h_2^{-2}(a,u_2,u_3).$ 
Therefore we obtain 
$\int_{\partial O_{\epsilon}} (h^{-2}_{1}-h^{-2}_{2})
(\partial_s \phi) \phi \nu_s =0.$ So (\ref{4.3}) becomes (\ref{4.6}) 
with $h_1^{-2}-h_2^{-2}>0$ in $O_{\epsilon}$ and $V_{k_1}-V_{k_2}>0$ in 
$O_{\epsilon}.$ So $\phi=0$ in $O_{\epsilon}$ and as in Step 1, by a unique 
continuation theorem, we obtain that $\phi=0$ in $\Omega_0$. 
Therefore we get a contradiction.

\end{document}